\documentclass[graybox]{svmult}

\usepackage{mathptmx}       
\usepackage{helvet}         
\usepackage{courier}        
\usepackage{type1cm}        
\usepackage{makeidx}         
\usepackage{graphicx}        
\usepackage{multicol}        
\usepackage[bottom]{footmisc}

\makeindex             

\newif\iffinal

\finaltrue

\usepackage{todonotes}
\usepackage{siunitx}
\usepackage[acronym]{glossaries}

\newacronym[longplural={interlevel set persistence hierarchies}]{ISPH}{ISPH}{interlevel set persistence hierarchy}

\usepackage{dsfont}

\usepackage{amsmath}
\usepackage{amssymb}

\DeclareMathOperator{\card}{card}
\DeclareMathOperator{\pers}{pers}
\DeclareMathOperator{\rank}{rank}
\DeclareMathOperator{\stab}{stab}

\newcommand{\domain}   {\ensuremath{\mathds{D}}}
\newcommand{\hierarchy}{\ensuremath{\mathcal{H}}}

\newcommand{\reaches}  {\ensuremath{\sim}}
\newcommand{\real}     {\ensuremath{\mathds{R}}}
\newcommand{\highest}  {\ensuremath{\hat{c}}}

\newcommand{\interlevelset}          [2]{\ensuremath{\mathcal{L}_{#1,#2}(f)}}
\newcommand{\levelset}               [1]{\ensuremath{\mathcal{L}_{#1}(f)}}
\newcommand{\sublevelset}            [1]{\ensuremath{\mathcal{L}_{#1}^{-}(f)}}

\newcommand{\landau}                 [1]{\ensuremath{\mathcal{O}\left(#1\right)}}

\usepackage{tikz}
\usepackage{pgfplots}

\pgfplotsset{compat=1.14}

\usepgfplotslibrary{colorbrewer}
\usepgfplotslibrary{fillbetween}

\definecolor{amber}     {RGB}{255,191,  0}
\definecolor{cardinal}  {RGB}{196, 30, 58}
\definecolor{yale}      {RGB}{ 70,130,180}

\pgfplotscreateplotcyclelist{color list}{%
  cardinal,
  yale,
  amber
}

\pgfplotsset{cycle list name=color list}

\usetikzlibrary{calc}
\usetikzlibrary{patterns}
\usetikzlibrary{shapes}

\tikzset{circle split part fill/.style  args={#1,#2}{%
  alias=tmp@name,
    postaction={%
      insert path={
       \pgfextra{%
         \pgfpointdiff{\pgfpointanchor{\pgf@node@name}{center}}%
                      {\pgfpointanchor{\pgf@node@name}{east}}%
         \pgfmathsetmacro\insiderad{\pgf@x}
           \fill[#1] (\pgf@node@name.base) ([xshift=-\pgflinewidth]\pgf@node@name.east) arc (  0:180:\insiderad-\pgflinewidth) -- cycle;
           \fill[#2] (\pgf@node@name.base) ([xshift= \pgflinewidth]\pgf@node@name.west) arc (180:360:\insiderad-\pgflinewidth) -- cycle;
}}}}}
\makeatother

\usepackage{subfigure}

\newcommand{\subfigureCaptionSkip}{\vspace{-10pt}}%

\usepackage{algorithm}
\usepackage{algorithmicx}
\usepackage{algpseudocode}

\makeatletter
  \algrenewcommand\ALG@beginalgorithmic{\small}
\makeatother

\floatname{algorithm}{Alg.}

\newcommand*\circled[1]{\tikz[baseline=(char.base)]{
    \node[shape=circle,fill=white,draw,inner sep=1pt] (char) {#1};}}
\newcommand*\circledT[1]{{\scriptsize\ensuremath{\protect\circled{#1}}}}

\newcommand{\comment}[1]{}
\usepackage[absolute]{textpos}

\begin{document}

\title*{Hierarchies and Ranks for Persistence Pairs}

\author{Bastian Rieck \and Filip Sadlo \and Heike Leitte}

\institute{%
Bastian Rieck \and Heike Leitte \at TU Kaiserslautern, \email{\{rieck, leitte\}@cs.uni-kl.de}%
\and Filip Sadlo \at Heidelberg University, \email{sadlo@uni-heidelberg.de}%
}

\maketitle

\abstract{%
  We develop a novel hierarchy for zero-dimensional persistence pairs, i.e., connected components,
  which is capable of capturing more fine-grained spatial relations between persistence pairs.
  Our work is motivated by a lack of spatial relationships between features in persistence diagrams,
  leading to a limited expressive power.
  We build upon a recently-introduced hierarchy of pairs in persistence diagrams that augments the
  pairing stored in persistence diagrams with information about \emph{which} components merge. Our
  proposed hierarchy captures differences in branching structure.
  Moreover, we show how to use our hierarchy to measure the spatial stability of a pairing and we
  define a rank function for persistence pairs and demonstrate different applications.
}

\iffinal
\else
  \begin{textblock*}{\paperwidth}[1, 0](\paperwidth,0.5cm)
  \scriptsize%
  Authors' copy. Please refer to \emph{Topological Methods in Data Analysis and Visualization
  V: Theory, Algorithms, and Applications} for the definitive version of this chapter.
  \end{textblock*}
\fi

\section{Introduction}

A wide range of application domains employ the concept of persistence, i.e., a measure of feature
robustness or scale. It is particularly effective when dealing with noisy data, permitting analysts
to distinguish between ``signal'' and ``noise''.
Being a purely topological approach, however, the information conferred by persistence does
not retain any spatial information.  While this is sometimes desirable, previous
work~\cite{Gerber10,Rieck14b,Zomorodian08} has shown that retaining at least a minimum of
geometrical information is often beneficial, as it increases the expressive power.
In this paper, we develop a hierarchy that relates points in a persistence diagram. Our hierarchy
makes exclusive use of topological properties of data, while still being able to distinguish between
geometrically distinct data. Moreover, the hierarchy is capable of measuring stability
properties of the pairing of critical points itself, yielding additional structural stability
information about input data.
We demonstrate the practicality of our method by means of several datasets. Additionally, we
compare it to a state-of-the-art hierarchy, point out the improvements over said hierarchy, and
demonstrate how our novel approach differs from related hierarchical concepts such as Reeb graphs.

\section{Related Work}

We refer the reader to Edelsbrunner and Harer~\cite{Edelsbrunner10} for
a detailed overview of persistence and related concepts.
There are several related approaches for creating a hierarchy of persistence information.
Doraiswamy et al.~\cite{Doraiswamy13} calculate a topological saliency of critical points in
a scalar field based on their spatial arrangement. Critical points with low persistence that are
isolated from other critical points have a higher saliency in this concept.
These calculations yield saliency curves for different smoothing radii. While these curves permit
a ranking of persistence pairs, they do not afford a description of their nesting behavior.
Consequently, in contrast to our approach, the saliency approach is incapable of distinguishing some
spatial rearrangements that leave persistence values and relative distances largely intact, such as
moving all peaks towards each other.
Bauer~\cite{Bauer11} developed what we refer to in this paper as the
regular persistence hierarchy.
It is fully combinatorial and merely requires
small changes of the pairing calculation of related critical points.
This hierarchy was successfully used in determining cancellation
sequences of critical points of surfaces.
However, as shown in this paper, this hierarchy cannot distinguish
between certain nesting relations.

In scalar field analysis, the calculation of graph structures such as the Reeb
graph~\cite{Doraiswamy09} or the contour tree~\cite{Carr03}, along with merge and split trees, has a long tradition.
These graphs relate critical points with each other, but do not permit the calculation of
hierarchies of persistence pairs. We will demonstrate this on a simple one-dimensional example in
this paper.
Recent work in this area is driven by the same motivation as our work: the merge tree,
for example, turns out to be more expressive with respect to spatial differences in the domain.
Thus, even if two scalar fields have the same critical pairs, their merge trees are capable of
retaining differences in sublevel set merging behavior.
This observation led to the development of distance measures for merge trees~\cite{Beketayev14},
Reeb graphs~\cite{Bauer14}, and extremum graphs~\cite{Narayanan15}.
Since the aforementioned tree structures tend to be unwieldy, Pascucci et
al.~\cite{Pascucci09} proposed a hierarchical decomposition, the branch decomposition. This
decomposition relates the different branches of a contour tree with each other.
Recently, Saikia et al.~\cite{Saikia14} used these graphs as a similarity measure for the structural
comparison of scalar data. While these works are close to our method in spirit, they rely on
a different type of structural information.

\section{Background and Notation}

In this paper, we assume that we are working with a domain~$\domain$ and a scalar function~$f\colon\domain\to\real$.
We make no assumptions about the connectivity of~$\domain$ or its dimension.
As for~$f$, we require it to have a finite number of critical points---a
condition that is always satisfied for real-world data---and
that the function values of those critical points are different---a
condition that may be satisfied by, e.g., symbolic
perturbation~\cite{Edelsbrunner90}.
Such scalar fields commonly occur in many different applications, and their features are often
described using scalar field topology.
This umbrella term refers to the analysis of how certain special sets---the level sets---of
the scalar function~$f$ change upon varying parameters.
More precisely, given a threshold~$c$, one typically distinguishes between, e.g., level set~$\levelset{c}$,
and sublevel set~$\sublevelset{c}$,
\begin{align}
  \levelset{c}      :=& \left\{ x \in \domain \mid f\left(x\right) = c \right\}\\
  \sublevelset{c}   :=& \left\{ x \in \domain \mid f\left(x\right) \leq c \right\}
\end{align}
In this paper, we also require the interlevel set~$\interlevelset{l}{u}$,
\begin{equation}
  \interlevelset{l}{u} := \sublevelset{u} \setminus \sublevelset{l} = \{ x \in \domain \mid l \leq
  f(x) \leq u \}.
\end{equation}
Interlevel sets are commonly used to describe the topology of
real-valued functions~\cite{Carlsson09a} or the robustness of homology
classes~\cite{Bendich10,Bendich13}.
Scalar field topology refers to the investigation of changes in connectivity in these sets. Such
changes are intricately connected to the critical points of~$f$ by means of Morse
theory~\cite{Milnor63}.
Focusing on the sublevel sets~(the case for superlevel sets can be solved by duality arguments), we find that (local) minima \emph{create} new connected components in the sublevel set,
while~(local) maxima---or saddles in higher dimensions---are responsible for merging two connected
components, thereby \emph{destroying} one of them.
Related creators and destroyers may thus be paired with each other~(using, e.g., the ``elder
rule''~\cite[p.\ 150]{Edelsbrunner10} that merges the connected component with
a higher---younger---function value into the one with a lower---older---function
value), which permits their use in various data structures. 

The persistence diagram is a generic data structure to represent such a pairing.
For every creator--destroyer pair, it contains a point in~$\real^2$ according to the corresponding
function values.
Persistence diagrams have many desirable stability properties~\cite{Cohen-Steiner07,Cohen-Steiner10}
and permit the calculation of different metrics. Unfortunately, they are sometimes too coarse to
describe both the topology \emph{and} geometry of a scalar function.
Given a point $(c,d)$ in a persistence diagram, where $c$ is the function value of the creator and
$d$ is the function value of the paired destroyer, the absolute difference $|d-c|$ is referred to as
the persistence~$\pers(c,d)$ of the pair.
Persistence permits a way to define whether certain pairs are more prominent than others. Roughly
speaking, the persistence of such a pair is the magnitude of the smallest perturbation that is able
to cancel it.

\section{Persistence Hierarchies}

The calculation of persistence always underlies the idea of a pairing, i.e., a way of relating
different parts of a function with each other.
Here, we shall only focus on zero-dimensional persistent homology, which
describes connected components in the sublevel sets of a function, and
the elder rule for pairing points.
Consequently, we have a relationship between local minima and local
maxima~(or saddles in higher dimensions) of a function. We leave the
treatment of other topological features for future work.
Moreover, we only cover the case of sublevel sets; superlevel set
calculations follow by analogy.

\subsection{Regular Persistence Hierarchy}

Bauer~\cite{Bauer11} observes that the process of merging two connected components permits the
definition of a natural hierarchy between persistence pairs.
More precisely, assume that we are given two connected components $\sigma$ and
$\sigma^\prime$, each created at a local minimum.
If $\sigma$ merges into $\sigma^\prime$ at, e.g., a local maximum, we consider $\sigma^\prime$ to be
the parent of $\sigma$.
This relation is a necessary but not sufficient condition for finding out which pairs of
critical points of a Morse function cannot be canceled without affecting other points.
We call this hierarchy, which is equivalent to a merge tree~\cite{Carr03}, the regular persistence
hierarchy. Each of its nodes corresponds to a creator--destroyer pair. The hierarchy forms
a directed acyclic graph, i.e., a tree. This is a consequence of the assumption that the function
values at critical points are unique.
Thus, whenever a merge of two connected components takes place, the ``younger'' component is uniquely
defined. Moreover, a critical point cannot both create and destroy a topological
feature, so there cannot be any cycles in the hierarchy.
The regular persistence hierarchy has a natural root that corresponds to the global minimum, as the
connected component corresponding to this value is never merged.
%

\runinhead{Example and limitations}
%
The regular persistence hierarchy cannot distinguish some connectivity
relations: for example, Fig.~\ref{fig:Example regular persistence
hierarchy} depicts the regular persistence hierarchies for two simple
functions. We can see that the hierarchy is equal for both functions
even though their connectivity behavior is different.
More precisely, in the red function, the two persistence pairs are connected via two different
branches of the function, i.e., it is impossible to reach both minima without traversing a third
minimum---the global one---as the threshold of the sublevel sets is increased.
This difference in connectivity results in a different stability of the pairing. A perturbation of
the critical points at \circledT{z} and \circledT{c} in the blue function, for example, is
capable of changing the complete pairing: if we move the points to $f(4)=1.9$ and $f(5)=0.9$,
respectively, the pairing of the critical point \circledT{a} at $x=1$ will change, as well. The same
perturbation has no effect on the red function, though.
A hierarchy of persistence pairs should account for these differences in connectivity.

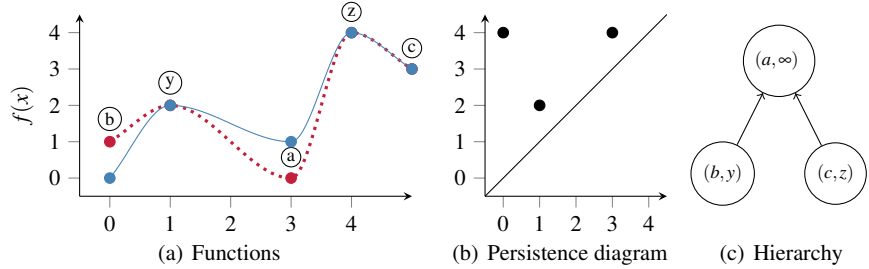
\begin{figure}[btp]
  \centering
  \subfigure[Functions]{%
    \iffinal
      \includegraphics{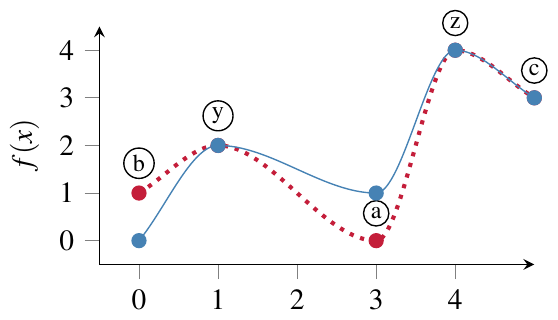}
    \else
      \begin{tikzpicture}
        \begin{axis}[%
          axis x line = bottom,
          axis y line = left,
          width       = 6cm,
          height      = 4cm,
          xtick       = {0,1,2,3,4},
          ytick       = {0,1,2,3,4},
          xmin        = -0.5,
          ymin        = -0.5,
          ymax        =  4.5,
          tick align  =  outside,
          ylabel      = $f(x)$,
        ]
          \addplot gnuplot [cardinal, no marks, raw gnuplot, dotted, very thick] {%
            set samples 1000;
            plot 'Data/Regular_persistence_hierarchy_f1.txt' smooth mcsplines;
          };
          \addplot[nodes near coords,
            every node near coord/.append style={xshift=0pt,yshift=1pt,anchor=south,font=\scriptsize},
            black, only marks, point meta=explicit symbolic] file {Data/Regular_persistence_hierarchy_f1.txt};
          \addplot[cardinal, only marks] file {Data/Regular_persistence_hierarchy_f1.txt};
          \addplot gnuplot [yale, no marks, raw gnuplot] {%
            set samples 1000;
            plot 'Data/Regular_persistence_hierarchy_f2.txt' smooth mcsplines;
          };
          \addplot[yale, only marks] file {Data/Regular_persistence_hierarchy_f2.txt};
        \end{axis}
      \end{tikzpicture}
    \fi
  }
  \hfill
  \subfigure[Persistence diagram]{%
    \iffinal
      \includegraphics{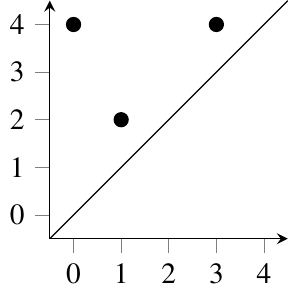}
    \else
      \begin{tikzpicture}
        \begin{axis}[%
          axis x line = bottom,
          axis y line = left,
          width       = 4cm,
          height      = 4cm,
          xmin        = -0.5,
          xmax        =  4.5,
          ymin        = -0.5,
          ymax        =  4.5,
          tick align  = outside,
          xtick       = {0,1,2,3,4},
          ytick       = {0,1,2,3,4},
        ]
          \addplot[only marks, black] coordinates {
            (0,4)
            (1,2)
            (3,4)
          };
          \addplot[domain=-0.5:4.5] {x};
        \end{axis}
      \end{tikzpicture}
    \fi
  }
  \hfill
  \subfigure[Hierarchy]{%
  \raisebox{0.18\height}{
        \iffinal
          \includegraphics{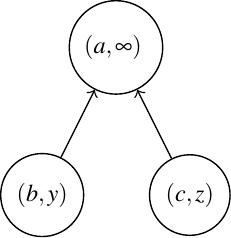}
        \else
          \begin{tikzpicture}[<-,every node/.style = {shape=circle, draw, align=center, font=\scriptsize}]

          \node { $(a,\infty)$ }
            child { node {$(b,y)$} }
            child { node {$(c,z)$} };
          \end{tikzpicture}
        \fi
      }
  }
  \subfigureCaptionSkip
  \caption{%
    Two functions with different connectivity but equal persistence diagrams. Both functions also
    share the same regular persistence hierarchy.
  }
  \label{fig:Example regular persistence hierarchy}
\end{figure}

\subsection{Interlevel Set Persistence Hierarchy~(ISPH)}

\begin{figure}[tbp]
  \centering
  \pgfdeclarelayer{background}
  \pgfsetlayers{background,main}
  \subfigure[\label{fig:Regular f1}]{%
    \iffinal
      \includegraphics{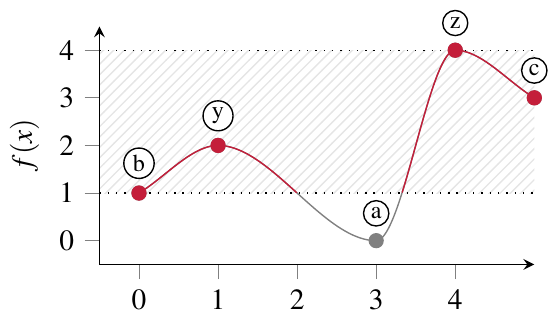}
    \else
      \begin{tikzpicture}
        \begin{axis}[%
          axis x line = bottom,
          axis y line = left,
          width       = 6cm,
          height      = 4cm,
          tick align  = outside,
          xtick       = {0,1,2,3,4},
          ytick       = {0,1,2,3,4},
          xmin        = -0.5,
          ymin        = -0.5,
          ymax        =  4.5,
          ylabel      = {$f(x)$},
        ]
          \addplot gnuplot [gray, no marks, raw gnuplot] {%
            set samples 1000;
            plot 'Data/Regular_persistence_hierarchy_f1.txt' smooth mcsplines;
          };
          \addplot[gray, only marks] file {Data/Regular_persistence_hierarchy_f1.txt};
          \addplot[nodes near coords,
            every node near coord/.append style={xshift=0pt,yshift=1pt,anchor=south,font=\scriptsize},
            black, only marks, point meta=explicit symbolic] file {Data/Regular_persistence_hierarchy_f1.txt};
          \addplot gnuplot [cardinal, no marks, raw gnuplot, restrict y to domain = 1:4] {%
            set samples 1000;
            plot 'Data/Regular_persistence_hierarchy_f1.txt' smooth mcsplines;
          };
          \addplot[cardinal, only marks, restrict y to domain = 1:4] file {Data/Regular_persistence_hierarchy_f1.txt};
          \draw[dotted,name path=Lower] (-0.5,1) -- (5,1);
          \draw[dotted,name path=Upper] (-0.5,4) -- (5,4);
          \begin{pgfonlayer}{background}
            \fill[pattern color=black!10, pattern=north east lines] (-0.5,1) rectangle (5,4);
          \end{pgfonlayer}
        \end{axis}
      \end{tikzpicture}
    \fi
  }
  \subfigure[\label{fig:Regular f2}]{%
    \iffinal
      \includegraphics{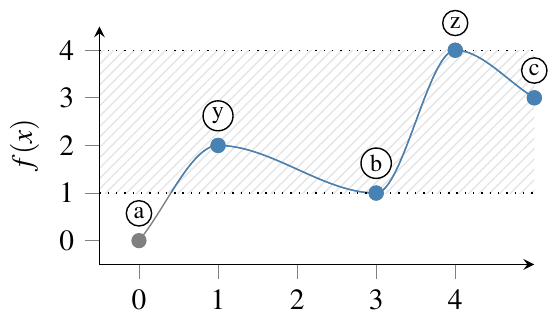}
    \else
      \begin{tikzpicture}
        \begin{axis}[%
          axis x line = bottom,
          axis y line = left,
          width       = 6cm,
          height      = 4cm,
          tick align  = outside,
          xtick       = {0,1,2,3,4},
          ytick       = {0,1,2,3,4},
          xmin        = -0.5,
          ymin        = -0.5,
          ymax        =  4.5,
          ylabel      = {$f(x)$},
        ]
          \addplot gnuplot [gray, no marks, raw gnuplot] {%
            set samples 1000;
            plot 'Data/Regular_persistence_hierarchy_f2.txt' smooth mcsplines;
          };
          \addplot[gray, only marks] file {Data/Regular_persistence_hierarchy_f2.txt};
          \addplot[nodes near coords,
            every node near coord/.append style={xshift=0pt,yshift=1pt,anchor=south,font=\scriptsize},
            black, only marks, point meta=explicit symbolic] file {Data/Regular_persistence_hierarchy_f2.txt};
          \addplot gnuplot [yale, no marks, raw gnuplot, restrict y to domain = 1:4] {%
            set samples 1000;
            plot 'Data/Regular_persistence_hierarchy_f2.txt' smooth mcsplines;
          };
          \addplot[yale, only marks, restrict y to domain = 1:4] file {Data/Regular_persistence_hierarchy_f2.txt};
          \draw[dotted,name path=Lower] (-0.5,1) -- (5,1);
          \draw[dotted,name path=Upper] (-0.5,4) -- (5,4);
          \begin{pgfonlayer}{background}
            \fill[pattern color=black!10, pattern=north east lines] (-0.5,1) rectangle (5,4);
          \end{pgfonlayer}

        \end{axis}
      \end{tikzpicture}
    \fi
  }
  \subfigureCaptionSkip
  \caption{%
    The merge phase of the \acrfull{ISPH} makes use of the connectivity of the
    interlevel set~(hatched lines): to connect the critical point pairs $(b,y)$ and $(c,z)$ in~\subref{fig:Regular
    f1}, a region belonging to a third critical point~\circledT{a} needs to be traversed.
    This is not the case for~\subref{fig:Regular f2}.
  }
  \label{fig:Improved persistence hierarchy merge}
\end{figure}
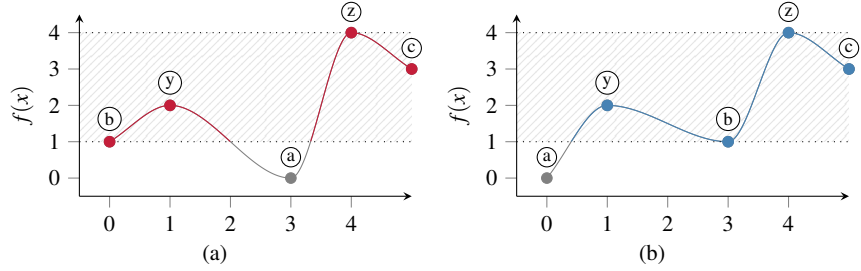

\begin{figure}[btp]
  \centering
  \hspace{2cm}
  \iffinal
    \includegraphics{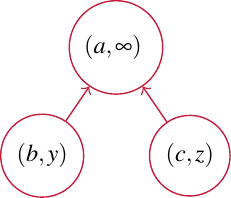}
  \else
    \begin{tikzpicture}[draw=cardinal,<-,every node/.style = {shape=circle, draw, align=center, font=\scriptsize}]
      \node { $(a,\infty)$ }
        child [level distance=1.1cm] { node {$(b,y)$} }
        child [level distance=1.1cm] { node {$(c,z)$} };
    \end{tikzpicture}
  \fi
  \hfill
  \raisebox{0.6\height}{
      \iffinal
        \includegraphics{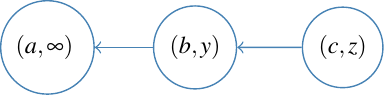}
      \else
        \begin{tikzpicture}[draw=yale,<-,every node/.style = {shape=circle, draw, align=center, font=\scriptsize},rotate=90]
          \node { $(a,\infty)$ }
          child { node {$(b,y)$}
            child { node {$(c,z)$} }
          };
        \end{tikzpicture}
      \fi
    }
  \hspace{2cm}
  \caption{%
    The \glspl{ISPH} for the example functions shown in Fig.~\ref{fig:Example
    regular persistence hierarchy}. In contrast to the regular persistence hierarchy, our hierarchy
    is capable of discriminating between the two functions. The hierarchy on the right
    has been rotated for layout reasons.
  }
  \label{fig:Example improved persistence hierarchy}
\end{figure}
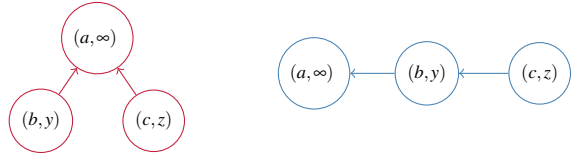

The example depicted in Fig.~\ref{fig:Example regular persistence hierarchy} demonstrated a lack of
discriminating information in the regular persistence hierarchy.
The key observation, illustrated as running example in Fig.~\ref{fig:Improved persistence hierarchy
merge}, is that not every merge of two connected components is topologically equal: a merge may
either result in a different branching structure of the hierarchy or it may keep the branches
intact.
Fig.~\ref{fig:Example improved persistence hierarchy} depicts our proposed \gls{ISPH}.
To measure these differences, we propose extending the traditional union--find data
structure that is used to detect the merges. Instead of merely storing the parent node of a given
connected component in the hierarchy, i.e., the generating critical point with lowest function
value, we also store \highest, the highest minimum---in terms of the function value---along this
branch. This will permit us to decide whether an additional branch needs to be introduced, as is
the case for function~\subref{fig:Regular f1} in Fig.~\ref{fig:Improved persistence hierarchy merge}.
If $\highest$ of a connected component is identical to the value of the parent, we call this
assignment trivial. We will use~$\highest$ interchangeably both for the critical point as well as
for its function value.
Subsequently, we distinguish between two types of merges: the first type of merge only extends
a branch in the hierarchy, while the second type of merge results in two branches that need to be
unified at a third critical point.
Fig.~\ref{fig:Improved persistence hierarchy merge} depicts the two cases for the example functions
shown in Fig.~\ref{fig:Example regular persistence hierarchy}.
We can see that for function~\subref{fig:Regular f1}, the pairs of critical points are connected by
an additional critical point~\circledT{a} only. Hence, two branches of the hierarchy merge at this point.
Function~\subref{fig:Regular f2}, by contrast, merely prolongs a branch in the hierarchy---both
pairs of critical points of the function are already connected without the inclusion of an
additional critical point.

To distinguish between these two cases, we check the stored highest critical points $\highest_l$ and
$\highest_r$ of the two connected components that merge at a local extremum.
Without loss of generality, we assume that $\highest_l$ belongs to the ``older'' branch and
$\highest_r$ belongs to the ``younger'' branch.
If both $\highest_l$ and $\highest_r$ are trivial, we merge their respective branches just as for
the regular persistence hierarchy.
Else, we have to check the induced connectivity to decide how the branches should be
connected. 
To this end, let $y_u$ refer to the value of the current critical point, i.e., the one at which the
two connected components merge. Furthermore, let $y_l$ refer to $\min(\highest_l, \highest_r)$, the
oldest of the two stored critical points.
We now calculate the interlevel set $\interlevelset{y_l}{y_u}$; see the colored parts in
Fig.~\ref{fig:Improved persistence hierarchy merge} for an example.
Following this, we check whether $\highest_l$ and $\highest_r$ are in the same
connected component with respect to $\interlevelset{y_l}{y_u}$. If so, the current branch can
be prolonged and $\highest_l$ is set to $\highest_r$. If not, two branches meet at the current
critical point and merge into one.
In Fig.~\ref{fig:Regular f1}, upon reaching \circledT{y}, we merge components \circledT{b} and
\circledT{a}, giving rise to the pair $(b,y)$. We have trivial critical points, i.e., $\highest_l
= \circledT{a}$ and $\highest_r = \circledT{b}$, so we add an edge between $(b,y)$ and $(a,\cdot)$
in the hierarchy; we do not yet know how $a$ will be paired, so we write ``$\cdot$'' for its
partner.
The next merge happens at \circledT{z}. We have $\highest_l = \circledT{b}$, $\highest_r
= \circledT{c}$, and $y_u = \circledT{z}$.  This gives rise to the interlevel set
$\interlevelset{b}{z}$. We now check whether $\highest_l$ and $\highest_r$ are connected in
$\interlevelset{b}{z}$. As this is not the case, we add an edge between $(c,z)$ and $(a,\cdot)$,
which is the parent of the older component, to the hierarchy.
In Fig.~\ref{fig:Regular f2}, by contrast, we have the same merges and the same interlevel set, but
\circledT{c} and \circledT{b} are \emph{connected} in $\interlevelset{b}{z}$, leading to the creation of an edge between $(c,z)$ and $(b,y)$.
The check with respect to the interlevel set connectivity is insufficient, however, for
higher-dimensional domains. Instead, we need to check whether a path in the neighborhood graph of
our data~(or, equivalently, an integral line) connecting the two critical points does not cross any
regions that are assigned to another critical point.
This presumes that we classify our data according to ascending or descending regions, which can be
easily integrated into standard persistent homology algorithms~\cite{Chazal13a}.

\runinhead{Algorithm and example}
%
The \gls{ISPH} requires only an additional data structure for storing information about
the critical points we encounter.
Moreover, we require checking the connectivity of the interlevel set---an operation that requires an
additional union--find data structure---and a way to calculate~(shortest) paths in the neighborhood
graph of our data. Alg.~\ref{alg:Extended persistence hierarchy} gives the pseudocode description of
our novel hierarchy based on sublevel sets.
Fig.~\ref{fig:Example improved persistence hierarchy} shows the \glspl{ISPH} for
the example functions in Fig.~\ref{fig:Example regular persistence hierarchy}, demonstrating
that our hierarchy represents differences in merging behavior in the sublevel sets.

\begin{algorithm}[tbp]
  \caption{Calculation of the \gls{ISPH}}
  \label{alg:Extended persistence hierarchy}
  \begin{algorithmic}[1]
    \newcommand{\UF}{\mathrm{U}}
    \Require A domain~$\domain$
    \Require A function~$f\colon\domain\to\real$
    \State $\UF \gets \emptyset$
    \State Sort the function values of~$f$ in ascending order
    \For{function value~$y$ of~$f$}
      \If{$y$ is a local minimum}
        \State Create a new connected component in $\UF$
        \State $\UF.\text{critical} \gets y$
      \ElsIf{$y$ is a local maximum or saddle}
        \State Use $\UF$ to merge the two connected components meeting at $y$
        \State Let $C'$ and $C$ be the two components meeting at $y$
        \If{both components have a trivial critical value}
          \State Create the edge $(C', C)$ in the hierarchy
        \Else
          \State Let $c'$ be the critical value of the older connected component
          \State Let $c$ be the critical value of the younger connected component
          \State $y_l \gets \min(c, c')$
          \State Create the interlevel set~$L := \interlevelset{y_l}{y}$
          \If{the shortest path connecting $c$, $c'$ in $L$ contains no other critical points}
            \State Create the edge $(c', y)$ in the hierarchy
          \Else
            \State Create the edge $(C', C)$ in the hierarchy~(as above)
          \EndIf
        \EndIf
      \Else%
        \State Use $\UF$ to add $y$ to the current connected component
      \EndIf
    \EndFor
  \end{algorithmic}

\end{algorithm}

\runinhead{Implementation}
%
We implemented our algorithm using \texttt{Aleph}\footnote{\url{https://github.com/Submanifold/Aleph}},
a library for exploring various uses of persistent homology. Our implementation of the \gls{ISPH} is
publicly available and supports processing structured grids~(using the VTK file format) as well as
one-dimensional functions.

\runinhead{Comparison with other tree-based concepts}

\begin{figure}[btp]
  \centering
  \subfigure[First function\label{sfig:Reeb example first function}]{%
    \iffinal
      \includegraphics{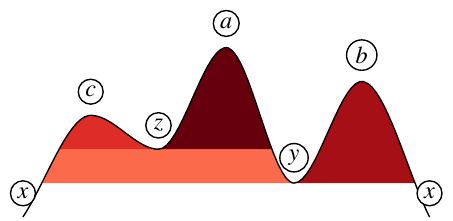}
    \else
      \begin{tikzpicture}
        \begin{axis}[%
          axis x line   = none,
          axis y line   = none,
          xmin          = -0.1,
          xmax          =  6.1,
          ymin          =  0.9,
          ymax          =  6.1,
          ytick         = \empty,
          xtick         = \empty,
          tick align    = outside,
          width         = 0.5\linewidth,
          unit vector ratio*= 2 1 1,
        ]
          \addplot gnuplot[black, no marks, raw gnuplot, name path=F] {%
            set samples 100;
            plot 'Data/Reeb_f1.txt' smooth mcsplines;
          };

          \addplot[%
            nodes near coords = {\scriptsize\circled{$\pgfplotspointmeta$}},
            point meta=explicit symbolic,
            only marks]
            table[header=false,meta index=2] {Data/Reeb_f1.txt};

          \path[name path=C1] (0,1) -- (6,1);
          \path[name path=C2] (0,2) -- (6,2);
          \path[name path=C3] (0,3) -- (6,3);

          \addplot[fill=none] fill between[of=C1 and F,
            split,
            every segment no 1/.style={fill=Reds-E},];

          \addplot[fill=none] fill between[of=C2 and F,
            split,
            every segment no 2/.style={fill=Reds-G},];

          \addplot[fill=none] fill between[of=C3 and F,
            split,
            every segment no 4/.style={fill=Reds-I},];

          \addplot[fill=none] fill between[of=C2 and F,
            split,
            every segment no 1/.style={fill=Reds-K},];

          \addplot[fill=none] fill between[of=C3 and F,
            split,
            every segment no 3/.style={fill=Reds-M},];
        \end{axis}
      \end{tikzpicture}
    \fi
  }
  \subfigure[Second function\label{sfig:Reeb example second function}]{%
    \iffinal
      \includegraphics{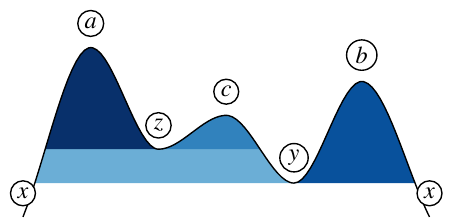}
    \else
      \begin{tikzpicture}
        \begin{axis}[%
          axis x line   = none,
          axis y line   = none,
          xmin          = -0.1,
          xmax          =  6.1,
          ymin          =  0.9,
          ymax          =  6.1,
          ytick         = \empty,
          xtick         = \empty,
          tick align    = outside,
          width         = 0.5\linewidth,
          unit vector ratio*= 2 1 1,
        ]
          \addplot gnuplot[black, no marks, raw gnuplot, name path=F] {%
            set samples 100;
            plot 'Data/Reeb_f2.txt' smooth mcsplines;
          };

          \addplot[%
            nodes near coords = {\scriptsize\circled{$\pgfplotspointmeta$}},
            point meta=explicit symbolic,
            only marks]
            table[header=false,meta index=2] {Data/Reeb_f2.txt};

          \path[name path=C1] (0,1) -- (6,1);
          \path[name path=C2] (0,2) -- (6,2);
          \path[name path=C3] (0,3) -- (6,3);

          \addplot[fill=none] fill between[of=C1 and F,
            split,
            every segment no 1/.style={fill=Blues-E},];

          \addplot[fill=none] fill between[of=C2 and F,
            split,
            every segment no 2/.style={fill=Blues-G},];

          \addplot[fill=none] fill between[of=C2 and F,
            split,
            every segment no 1/.style={fill=Blues-K},];

          \addplot[fill=none] fill between[of=C3 and F,
            split,
            every segment no 3/.style={fill=Blues-I},];

          \addplot[fill=none] fill between[of=C3 and F,
            split,
            every segment no 4/.style={fill=Blues-M},];
        \end{axis}
      \end{tikzpicture}
    \fi
  }
  \subfigure[Reeb graph\label{sfig:Reeb graph}]{%
    \iffinal
      \includegraphics{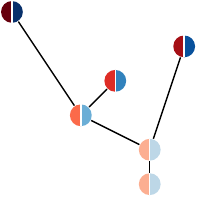}
    \else
      \begin{tikzpicture}[%
        scale = 0.35,
        every node/.style = {shape=circle split, draw=white, align=center, scale=.5, rotate=90},
      ]

        \node[circle split part fill={Reds-E,Blues-E}] (C1) at ( 0,1) {};
        \node[circle split part fill={Reds-E,Blues-E}] (C2) at ( 0,2) {};
        \node[circle split part fill={Reds-G,Blues-G}] (C3) at (-2,3) {};
        \node[circle split part fill={Reds-I,Blues-I}] (C4) at (-1,4) {};
        \node[circle split part fill={Reds-K,Blues-K}] (C5) at ( 1,5) {};
        \node[circle split part fill={Reds-M,Blues-M}] (C6) at (-4,6) {};

        \draw (C1) -- (C2);
        \draw (C2) -- (C5);
        \draw (C2) -- (C3);
        \draw (C3) -- (C4);
        \draw (C3) -- (C6);
      \end{tikzpicture}
    \fi
  }
  \subfigure[\Acrfullpl{ISPH}\label{sfig:Reeb example extended hierarchies}]{
    \iffinal
      \includegraphics{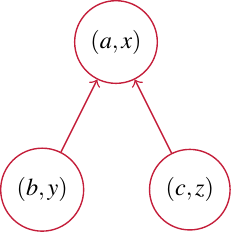}
    \else
      \begin{tikzpicture}[%
        <-,
        draw              = {cardinal},
        every node/.style = {shape=circle, draw, align=center, font=\scriptsize, minimum size=2.0em},
      ]
        \node {$(a,x)$}
          child {
            node {$(b,y)$}
          }
          child {
            node {$(c,z)$}
          };
      \end{tikzpicture}
    \fi
    \quad
    \raisebox{0.75\height}{
      \iffinal
        \includegraphics{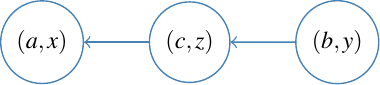}
      \else
        \begin{tikzpicture}[%
          <-,
          draw              = {yale},
          every node/.style = {shape=circle, draw, align=center, font=\scriptsize, minimum size=2.0em},
          rotate            = 90,
        ]
          \node {$(a,x)$}
            child {
              node {$(c,z)$}
              child {
                node {$(b,y)$}
              }
            };
        \end{tikzpicture}
      \fi
    }
  }
  \subfigureCaptionSkip
  \caption{%
    Comparison with Reeb graphs. The two functions yield the same Reeb graph~(or, equivalently, the merge tree of their superlevel sets), while our hierarchy is capable of telling them apart.
  }
  \label{fig:Reeb graph comparison}
\end{figure}
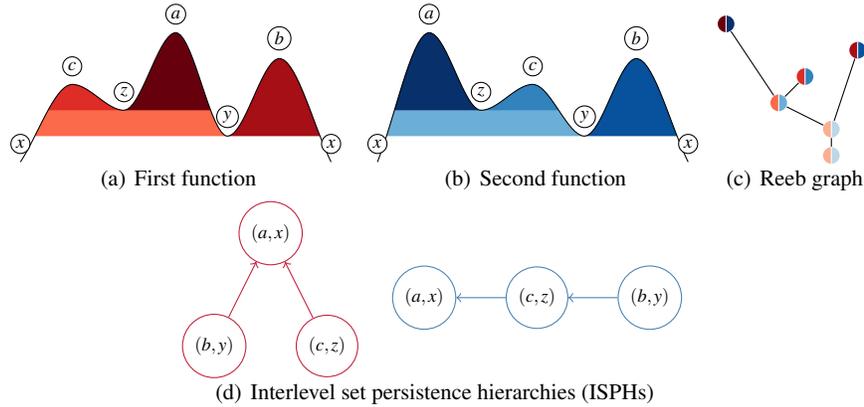

The \gls{ISPH} is capable of preserving more information than merge trees,
split trees, and Reeb graphs. As a simple example, consider the functions in Fig.~\ref{fig:Reeb
graph comparison}. Both functions carry the same sublevel/superlevel set information; their persistence
diagrams and regular persistence hierarchies coincide. Their Reeb graphs, shown in
Fig.~\ref{sfig:Reeb graph} with nodes whose colors indicate the corresponding contour of the
function, are also equal.
Of course, this does not imply that Reeb graphs~(or merge trees) are generally unsuitable.  During
our experiments, we encountered numerous functions in which Reeb graphs~(or merge trees) are able to
detect differences in functions with equal persistence diagrams. At the same time, the \gls{ISPH} was able to detect differences in these cases as well.
Fig.~\ref{sfig:Reeb example extended hierarchies} shows the \glspl{ISPH} of the
functions in Fig.~\ref{sfig:Reeb example first function} and Fig.~\ref{sfig:Reeb example second
function}.

The preceding example proves that the \gls{ISPH} is unrelated to existing decompositions: since the
Reeb graphs~(and the merge trees) of the two functions are equal but the \glspl{ISPH} differ, it is
not possible to derive the \gls{ISPH} from, e.g., the branch decomposition tree~\cite{Pascucci09} or the
extended branch decomposition graph~\cite{Saikia14}.

\runinhead{Robustness}
%
When adding noise to a function, topological hierarchies such as the merge~(split) tree and the Reeb
graph are known to contain numerous short branches, which make identifying important features and
comparing different trees more difficult~\cite{Saikia14}. By contrast, our novel \gls{ISPH} only contains as many points as there are \emph{pairs} in the persistence
diagram. Moreover, low-persistence pairs do not result in too much clutter because they tend to only
create short branches. In that sense, our hierarchy performs similarly well as the extended branch
decomposition graph by Saikia et al.~\cite{Saikia14}.

\subsubsection{Calculating Ranks}

Since the \gls{ISPH} is a DAG, we can define the rank of a topological
feature: given two vertices $u$ and $v$ in the hierarchy~$\hierarchy$, we write $u \reaches v$ if
there is a directed path connecting $u$ and $v$. The rank of a vertex $u$ in the hierarchy is then
calculated as the number of vertices that are reachable from it, i.e.,
\begin{equation}
  \rank(u) := \card\left\{ v \in \hierarchy \mid u \reaches v \right\},
\end{equation}
with $\rank(\cdot) \in \mathds{N}$. The minimum of the rank function is obtained for the last
connected component to be destroyed, i.e., the one that merges last with another component.
The rank can be easily calculated using a depth-first traversal of the tree.
It may be visualized as additional information within a persistence diagram, thereby
permitting datasets with similar persistence diagrams but different hierarchies to be distinguished
from each other without showing the hierarchy itself.
It is also invariant with respect to scaling of the function values in the data.
Similar concepts, such as the rank invariant~\cite{Carlsson09} in multidimensional persistence, only
use existing information from the persistence diagram, whereas our \gls{ISPH} goes beyond the
persistence diagram by including more structural information about critical points.

\subsubsection{Stability Measure}
\label{sec:Stability Measure}

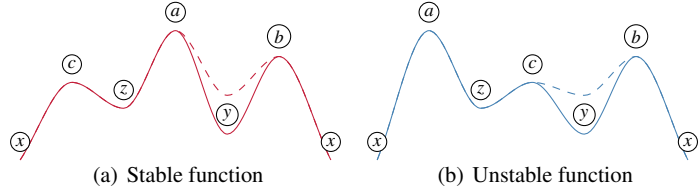
\begin{figure}[tbp]
  \centering
  \subfigure[Stable function]{%
    \iffinal
      \includegraphics{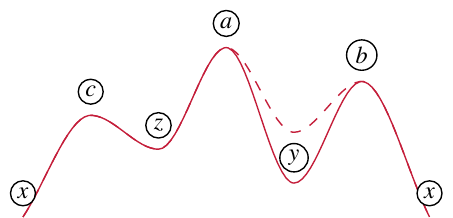}
    \else
      \begin{tikzpicture}
        \begin{axis}[%
          axis x line   = none,
          axis y line   = none,
          xmin          = -0.1,
          xmax          =  6.1,
          ymin          =  0.9,
          ymax          =  6.1,
          ytick         = \empty,
          xtick         = \empty,
          tick align    = outside,
          width         = 0.5\linewidth,
          unit vector ratio*= 2 1 1,
        ]
          \addplot gnuplot[cardinal, no marks, raw gnuplot, name path=F] {%
            set samples 100;
            plot 'Data/Reeb_f1.txt' smooth mcsplines;
          };
          \addplot gnuplot[cardinal, dashed, no marks, raw gnuplot, name path=F] {%
            set samples 100;
            plot 'Data/Reeb_f1_perturbed.txt' smooth mcsplines;
          };
          \addplot[%
            nodes near coords = {\scriptsize\circled{$\pgfplotspointmeta$}},
            point meta=explicit symbolic,
            only marks]
            table[header=false,meta index=2] {Data/Reeb_f1.txt};
        \end{axis}
      \end{tikzpicture}
    \fi
  }
  \subfigure[Unstable function]{%
    \iffinal
      \includegraphics{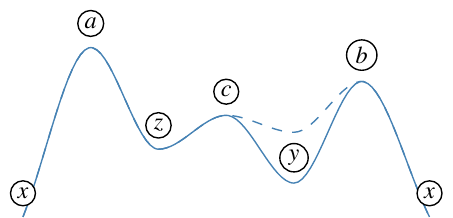}
    \else
      \begin{tikzpicture}
        \begin{axis}[%
          axis x line   = none,
          axis y line   = none,
          xmin          = -0.1,
          xmax          =  6.1,
          ymin          =  0.9,
          ymax          =  6.1,
          ytick         = \empty,
          xtick         = \empty,
          tick align    = outside,
          width         = 0.5\linewidth,
          unit vector ratio*= 2 1 1,
        ]
          \addplot gnuplot[yale, no marks, raw gnuplot, name path=F] {%
            set samples 100;
            plot 'Data/Reeb_f2.txt' smooth mcsplines;
          };

          \addplot gnuplot[yale, dashed, no marks, raw gnuplot, name path=F] {%
            set samples 100;
            plot 'Data/Reeb_f2_perturbed.txt' smooth mcsplines;
          };
          \addplot[%
            nodes near coords = {\scriptsize\circled{$\pgfplotspointmeta$}},
            point meta=explicit symbolic,
            only marks]
            table[header=false,meta index=2] {Data/Reeb_f2.txt};
        \end{axis}
      \end{tikzpicture}
    \fi
  }
  \subfigureCaptionSkip
  \caption{%
    Stable and unstable function whose superlevel sets yield the same persistence diagram. Since
    the \gls{ISPH} is capable of distinguishing between the two cases, it assesses their
    respective stability differently.
  }
  \label{fig:Stability illustration}
\end{figure}

Our \gls{ISPH} permits assessing the stability of the \emph{location} of critical points in the
pairing.
This issue with persistence diagrams was already pointed out by Bendich and
Bubenik~\cite{Bendich15}, who demonstrated that small changes in the critical values of
a function---while not drastically changing the persistence diagram itself---may still change the
points that are responsible for creating a certain topological feature.
The \gls{ISPH} contains information that may be used to assess the stability of
the creators of topological features.
To make this more precise, consider the example in Fig.~\ref{fig:Stability illustration}. Upon
traversing the superlevel sets of both functions, they exhibit the same persistence pairs, namely
$(a,z)$, $(b,y)$, and $(c,z)$. Refer to Fig.~\ref{sfig:Reeb example extended hierarchies} for the
corresponding \glspl{ISPH}.
If we perturb the critical point \circledT{y} for both
functions~(indicated by a dashed line), we still get the pairs $(b,y)$ and $(c,z)$ for the stable
function. For the unstable function, however, the perturbation results in the pairs $(b,z)$ and
$(c,y)$. In this sense, their location is less stable.

We thus define a stability measure for each critical point based on the hierarchy.
First, we need to quantify the stability of an \emph{edge}. Let $e := \{ (\sigma, \tau), (\sigma',\tau') \}$ be
an edge in the \gls{ISPH}. We define the stability of $e$ to be
\begin{equation}
  \stab(e) := \max\left\{|f(\sigma)- f(\sigma')|, |f(\tau) - f(\tau')| \right\},
\end{equation}
which is the minimum amount of perturbation that is required to change the hierarchy in the sense
described above. This quantity is also equal to the $\mathrm{L}_\infty$-distance between two points
in a persistence diagram.
We may now extend the stability measure to individual vertices by assigning each vertex~$v$ the
minimum stability value of its outgoing edges, i.e.,
\begin{equation}
  \stab(v) := \min\big\{\min\{ \stab(e) \mid e = (v,w) \in \hierarchy \}, \pers(v) \big\},
\end{equation}
where $w$ ranges over all direct children of $v$. Taking the second minimum ensures that we use the
persistence of $v$ if $v$ is a leaf node.

\subsubsection{Dissimilarity Measure}
\label{sec:Dissimilarity measure}

Since the \gls{ISPH} is a directed tree, a straightforward dissimilarity measure
is given by tree edit distance~\cite{Bille05} algorithms. These algorithms attempt to transform two
trees into each other by three elementary operations: relabeling a given node, deleting an existing
node, and inserting a new node.
Given two nodes with corresponding minima--maxima~$(c_1,d_1)$ and~$(c_2,d_2)$, respectively, we
define the cost for \emph{relabeling} to be
\begin{equation}
  \mathrm{cost}_1 = \max\big( |c_1-c_2|, |d_1-d_2| \big),
\end{equation}
i.e., the $\mathrm{L}_\infty$-distance between the two points. Similarly, we define the cost for
\emph{deleting} or \emph{inserting} a node somewhere else in the hierarchy to be
\begin{equation}
  \mathrm{cost}_2 = \pers(c,d) = |d-c|,
\end{equation}
i.e., the persistence of the pair. The choice of these costs is ``natural'' in the sense that they are
also used, e.g., when calculating the bottleneck distance~\cite{Cohen-Steiner07} between persistence
diagrams.

\runinhead{Complexity}
%
The tree edit distance has the advantage of being efficiently-solvable via standard dynamic
programming techniques. It is thus scalable with respect to the size of the hierarchy---more so than
the Wasserstein or bottleneck distances~\cite{Maria14, Tierny18} that are commonly used for comparing persistence diagrams:
we have a worst case complexity of~$\landau{n^2 m \log m}$, where $n$ is the number of pairs in the
smaller hierarchy, and $m$ is the number of pairs in the larger hierarchy. By contrast, the
Wasserstein distance has a time complexity of~$\landau{n^3}$~\cite[p.\ 196]{Edelsbrunner10}, and we observed large differences in
runtime behavior~(see Sec.~\ref{sec:Results}).

\section{Results}
\label{sec:Results}

We exemplify some usage scenarios of the \gls{ISPH} by means of several
synthetic and non-synthetic datasets.

\subsection{Synthetic Data}
%
\begin{figure}[tbp]
  \centering
  \subfigure[\label{fig:Three peaks stable}]{%
    \raisebox{0.5cm}{%
      \includegraphics[height=1.00cm]{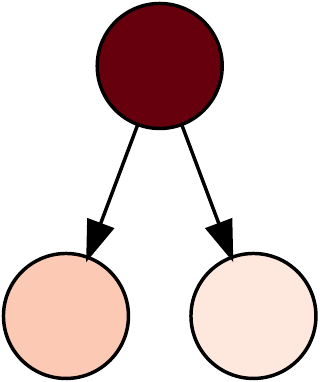}
    }
    \includegraphics[height=2.00cm]{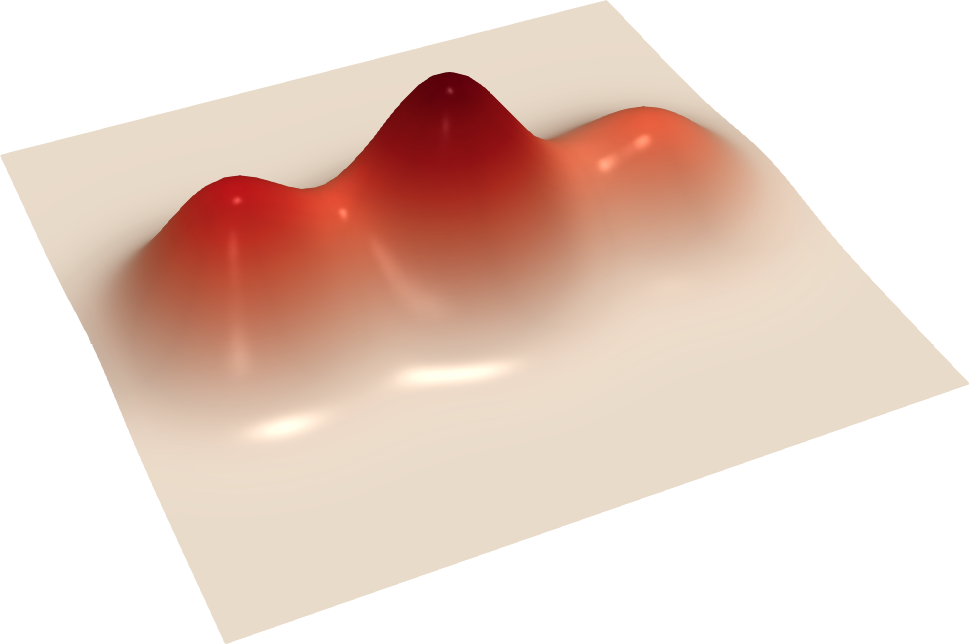}
  }
  \subfigure[\label{fig:Three peaks ridge}]{%
    \includegraphics[height=2cm]{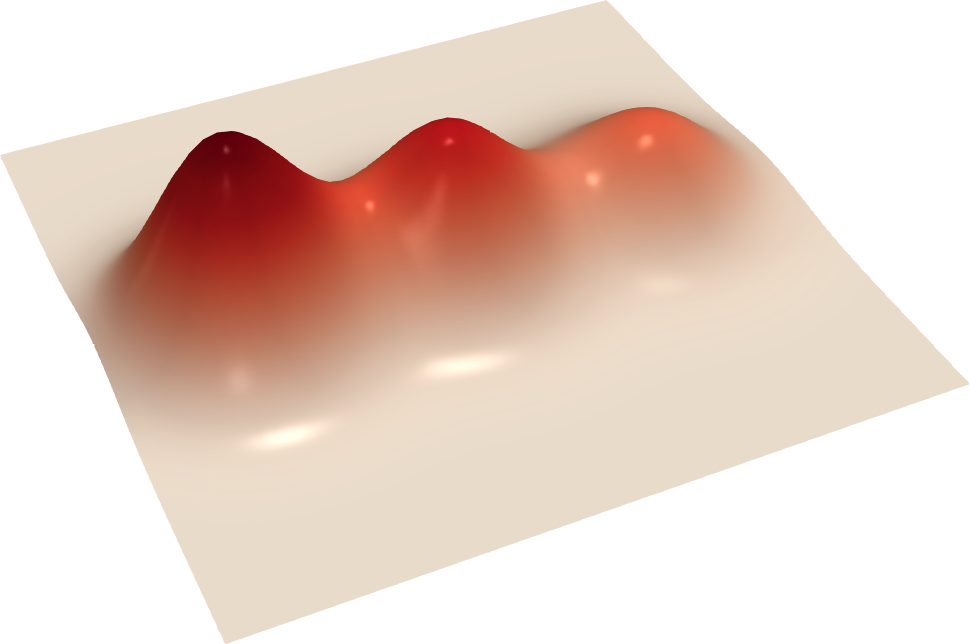}
    \raisebox{0.25cm}{%
      \includegraphics[height=1.5cm]{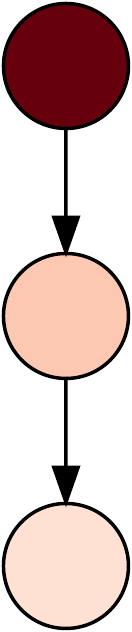}
    }
  }
  \subfigureCaptionSkip
  \caption{%
    In contrast to previous approaches, our hierarchy is capable of distinguishing between two peaks
    that are connected via a third one that is higher~\subref{fig:Three peaks stable} and a ``ridge''
    of peaks~\subref{fig:Three peaks ridge}.
  }
  \label{fig:Three peaks}
\end{figure}
\begin{figure}[btp]
  \centering
  \subfigure[Persistence diagrams for Fig.~\ref{fig:Three peaks}\label{fig:Three peaks persistence diagrams}]{%
    \iffinal
      \includegraphics{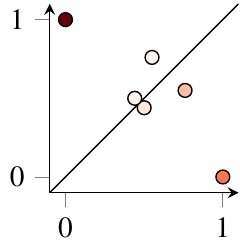}%
      \includegraphics{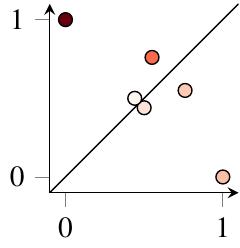}%
    \else
      \begin{tikzpicture}
        \begin{axis}[%
          scatter/use mapped color = {%
            draw = black!,
            fill = mapped color
          },
          height = 3.5cm,
          width  = 3.5cm,
          xmin =  -0.10,
          xmax =   1.10,
          ymin =  -0.10,
          ymax =   1.10,
          axis x line = bottom,
          axis y line = left,
          tick align  = outside,
          xtick  = { 0, 1},
          ytick  = { 0, 1},
          point meta min = 0.00,
          point meta max = 1.00,
          colormap/Reds-9,
        ]
          \addplot[only marks, scatter, point meta=explicit] file {Data/3_persistence_diagram_stability.txt};
          \addplot[only marks, scatter, point meta=explicit] file {Data/3_persistence_diagram_ranks.txt};
          \addplot[black, domain = {-0.10:1.10}] {x};
        \end{axis}
      \end{tikzpicture}
      \begin{tikzpicture}
        \begin{axis}[%
          scatter/use mapped color = {%
            draw = black!,
            fill = mapped color
          },
          height = 3.5cm,
          width  = 3.5cm,
          xmin = -0.10,
          xmax =  1.10,
          ymin = -0.10,
          ymax =  1.10,
          axis x line = bottom,
          axis y line = left,
          tick align  = outside,
          xtick  = { 0, 1},
          ytick  = { 0, 1},
          point meta min = 0.00,
          point meta max = 1.00,
          colormap/Reds-9,
        ]
          \addplot[only marks, scatter, point meta=explicit] file {Data/3_ridge_persistence_diagram_stability.txt};
          \addplot[only marks, scatter, point meta=explicit] file {Data/3_ridge_persistence_diagram_ranks.txt};
          \addplot[black, domain = {-0.10:1.10}] {x};
        \end{axis}
      \end{tikzpicture}
    \fi
  }
  \subfigure[Persistence diagrams for Fig.~\ref{fig:Peaks and craters}\label{fig:Peaks and craters persistence diagrams}]{%
    \iffinal
      \includegraphics{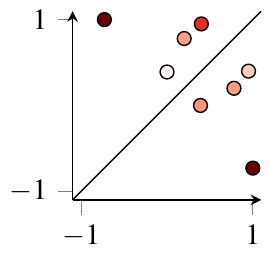}%
      \includegraphics{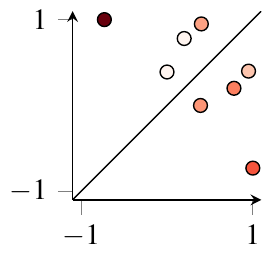}%
    \else
      \begin{tikzpicture}
        \begin{axis}[%
          scatter/use mapped color = {%
            draw = black!,
            fill = mapped color
          },
          height = 3.5cm,
          width  = 3.5cm,
          xmin = -1.10,
          xmax =  1.10,
          ymin = -1.10,
          ymax =  1.10,
          axis x line = bottom,
          axis y line = left,
          tick align  = outside,
          xtick  = { -1.0, 1},
          ytick  = { -1.0, 1},
          point meta min = 0.00,
          point meta max = 1.15,
          colormap/Reds-9,
        ]
          \addplot[only marks, scatter, point meta=explicit] file {Data/Peaks_craters_01_persistence_diagram_stability.txt};
          \addplot[only marks, scatter, point meta=explicit] file {Data/Peaks_craters_01_persistence_diagram_ranks.txt};
          \addplot[black, domain = {-1.10:1.10}] {x};
        \end{axis}
      \end{tikzpicture}
      \begin{tikzpicture}
        \begin{axis}[%
          scatter/use mapped color = {%
            draw = black!,
            fill = mapped color
          },
          height = 3.5cm,
          width  = 3.5cm,
          xmin = -1.10,
          xmax =  1.10,
          ymin = -1.10,
          ymax =  1.10,
          axis x line = bottom,
          axis y line = left,
          tick align  = outside,
          xtick  = { -1.0, 1},
          ytick  = { -1.0, 1},
          point meta min = 0.00,
          point meta max = 1.15,
          colormap/Reds-9,
        ]
          \addplot[only marks, scatter, point meta=explicit] file {Data/Peaks_craters_02_persistence_diagram_stability.txt};
          \addplot[only marks, scatter, point meta=explicit] file {Data/Peaks_craters_02_persistence_diagram_ranks.txt};
          \addplot[black, domain = {-1.10:1.10}] {x};
        \end{axis}
      \end{tikzpicture}
    \fi
  }
  \subfigureCaptionSkip
  \caption{%
    Combined persistence diagrams showing the ranks~(upper part) and the stability values~(lower
    part) as colors. Both carry sufficient information to distinguish datasets.
  }
  \label{fig:Synthetic data persistence diagrams}
\end{figure}
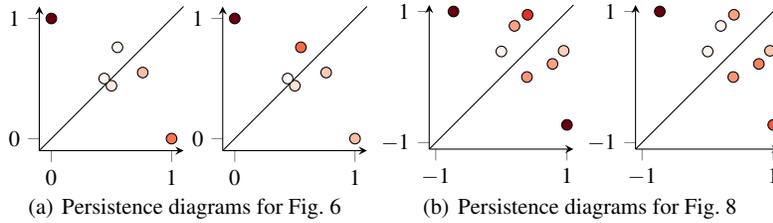

%
We created two synthetic datasets on a grid with 5,000 cells. Processing each dataset takes
approximately \SI{4.5}{\second}---regardless of whether we calculate the regular persistence hierarchy or the \gls{ISPH}. Our current implementation leaves lots of room for performance improvements,
though. 
Fig.~\ref{fig:Three peaks} shows the data together with the resulting hierarchies for the
two-dimensional equivalent to the data shown in Fig.~\ref{fig:Stability illustration}. We use
a standard color map that uses red for large values and white for low values.
Notice that the persistence diagrams of both datasets are equal, as well as their regular
persistence hierarchies~(which we do not show here).
Fig.~\ref{fig:Three peaks persistence diagrams} depicts the persistence diagrams of the two data
sets, colored by their stability values and the ranks~(mirrored along the diagonal).
Even this summary information gleaned from the \gls{ISPH} is capable of yielding
information to distinguish different datasets.

\begin{figure}
  \centering
  \subfigure[\label{fig:Peaks and craters 01}]{%
    \raisebox{0.5cm}{%
      \includegraphics[height=2cm]{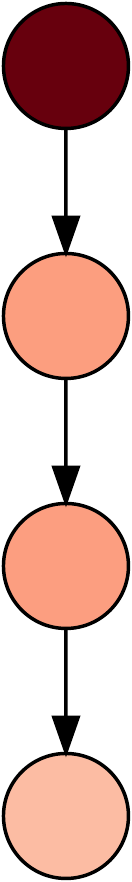}
    }
    \includegraphics[height=3cm]{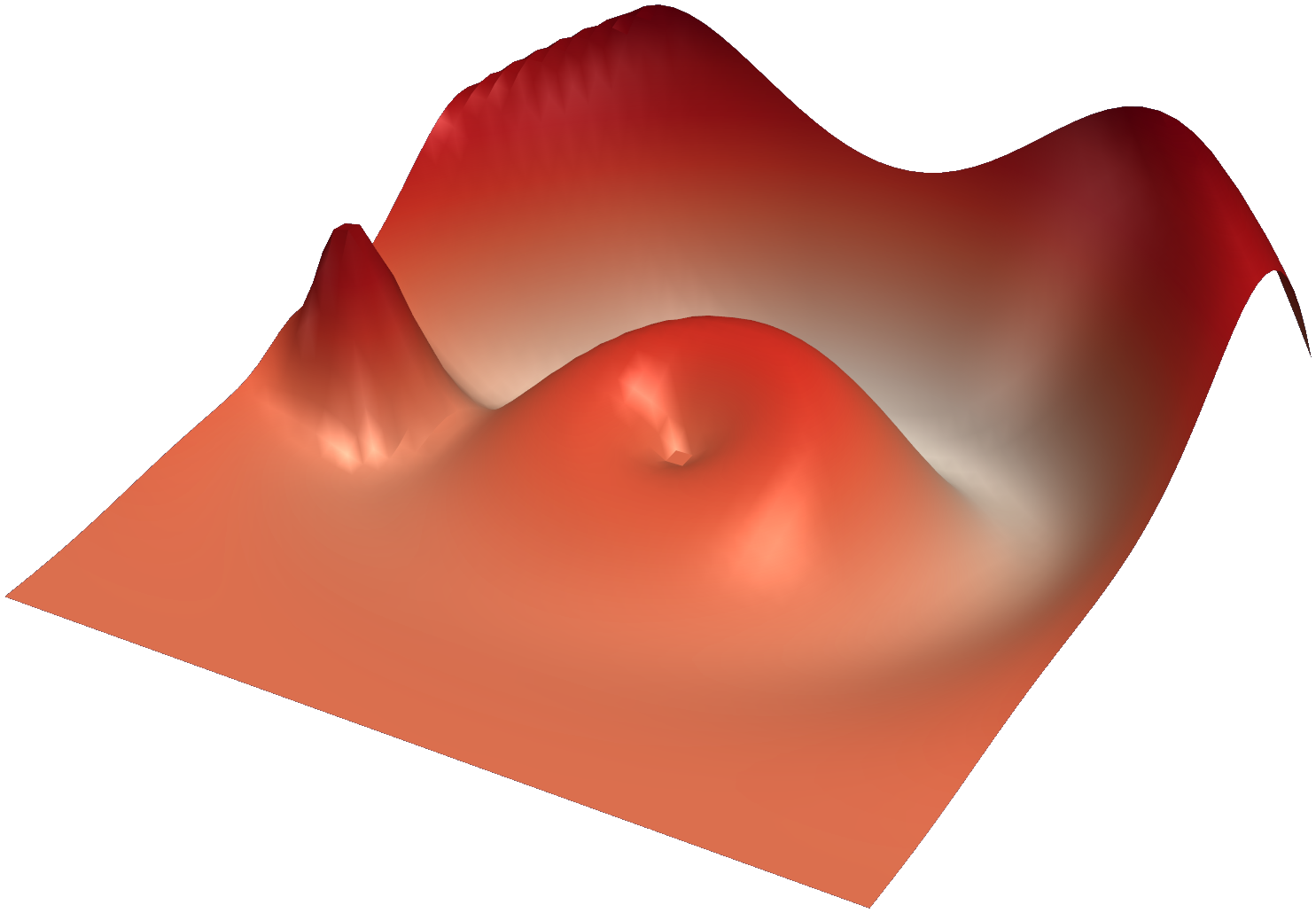}
  }
  \subfigure[\label{fig:Peaks and craters 02}]{%
    \includegraphics[height=3.0cm]{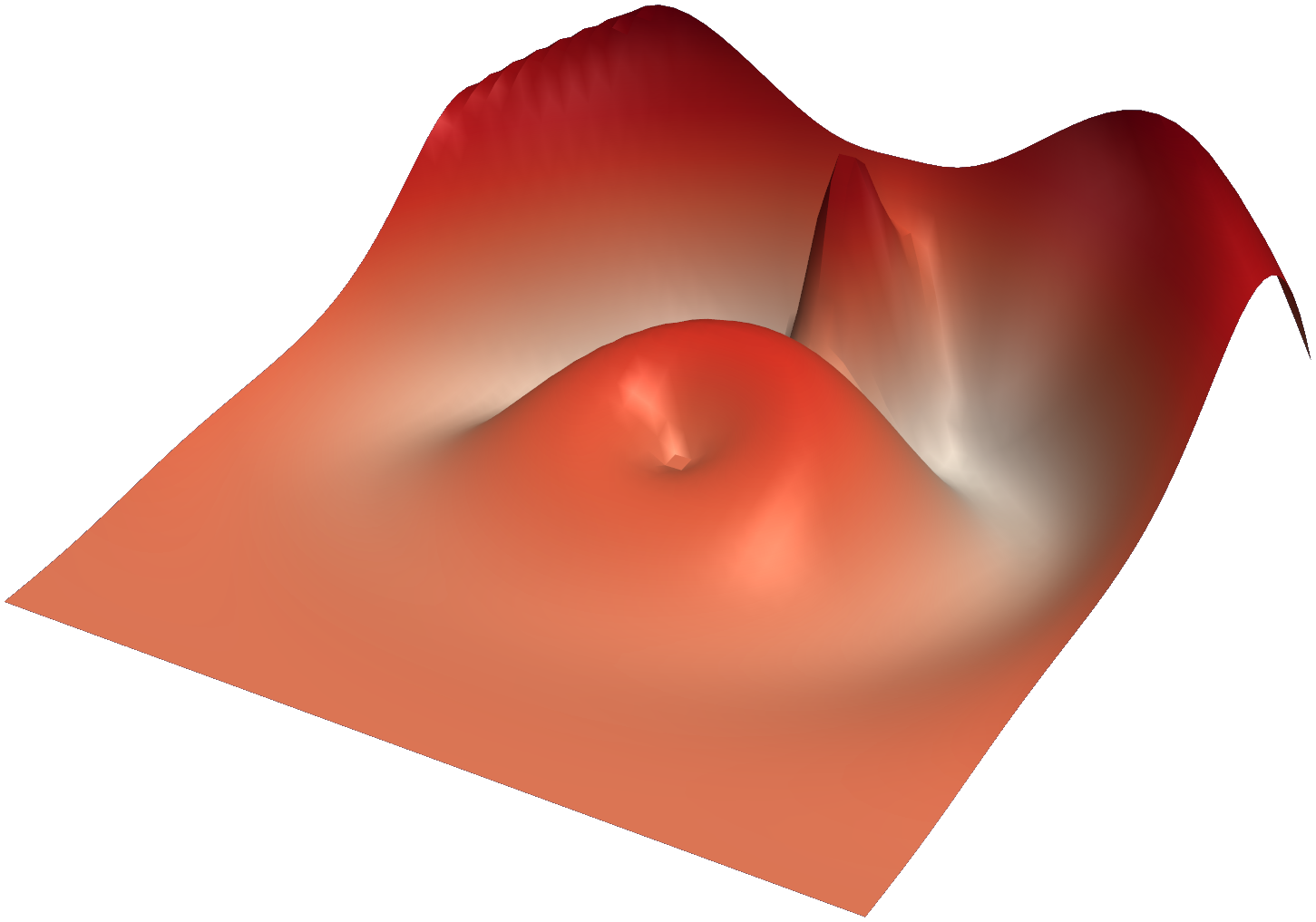}
    \raisebox{0.5cm}{%
      \includegraphics[height=1.5cm]{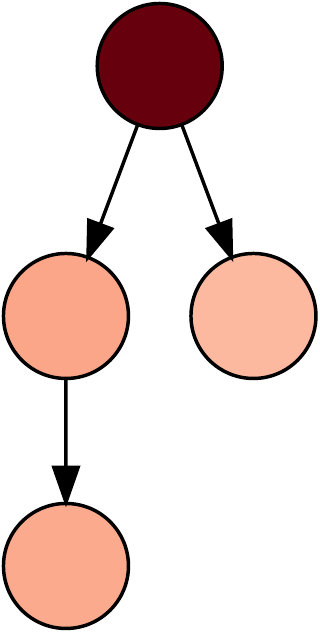}
    }
  }
  \subfigureCaptionSkip
  \caption{%
    Depending on the position of a single peak that moves, our hierarchies are different, because the
    merging behavior of critical points has been changed.
  }
  \label{fig:Peaks and craters}
\end{figure}

Fig.~\ref{fig:Peaks and craters} depicts a more complicated example, mixing peaks and craters. The
situation in Fig.~\ref{fig:Peaks and craters 02} is less stable because a perturbation of the peak
is capable of changing the complete pairing. Since this peak is connected differently in
Fig.~\ref{fig:Peaks and craters 01}, it does not decrease the stability of the global maximum. Note that the location of the peak is allowed to move; as long as it stays on the ridge, as depicted in Fig.~\ref{fig:Peaks and craters 02}, the \gls{ISPH} will not change. Likewise, as long as the peak moves along the plateau in the foreground, the \gls{ISPH} will remain the same as depicted in Fig.~\ref{fig:Peaks and craters 01}.

\subsection{Climate Data}
%
\begin{figure}[tbp]
  \centering
  \subfigure[$t=0$]{%
    \includegraphics[height=2.5cm]{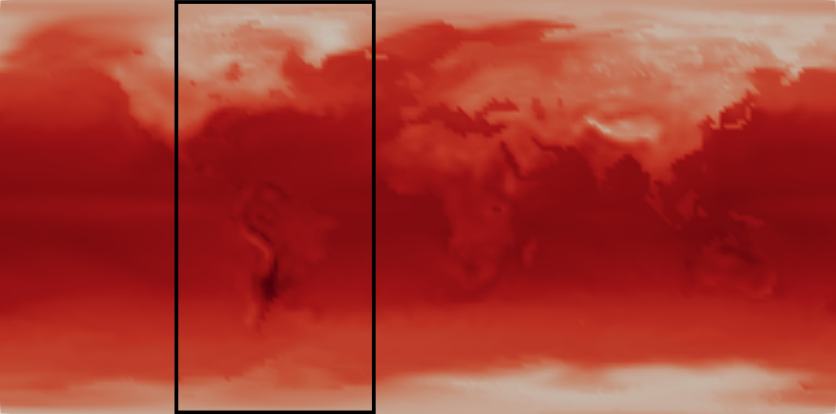}
  }
  \subfigure[$t=1$]{%
    \includegraphics[height=2.5cm]{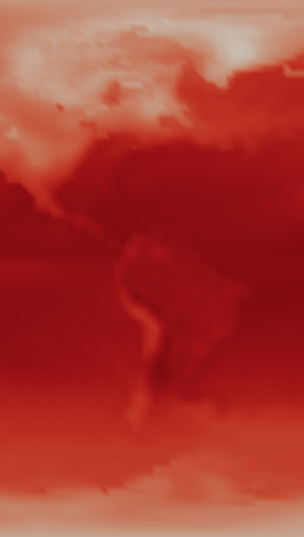}
  }
  \subfigure[$t=2$]{%
    \includegraphics[height=2.5cm]{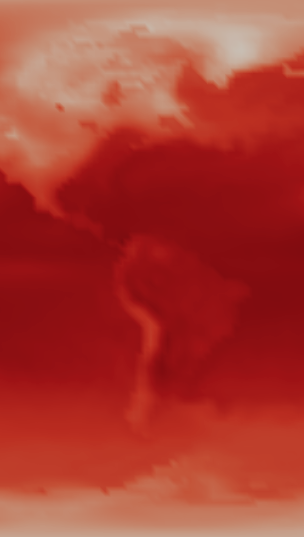}
  }
  \subfigure[$t=3$]{%
    \includegraphics[height=2.5cm]{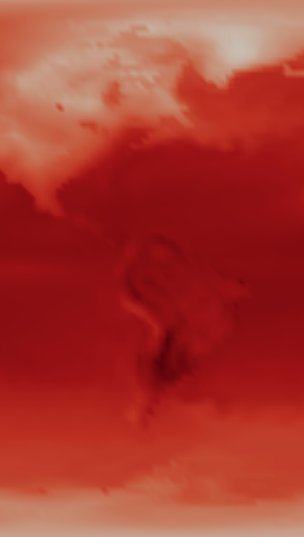}
  }
  \subfigure[$t=4$]{%
    \includegraphics[height=2.5cm]{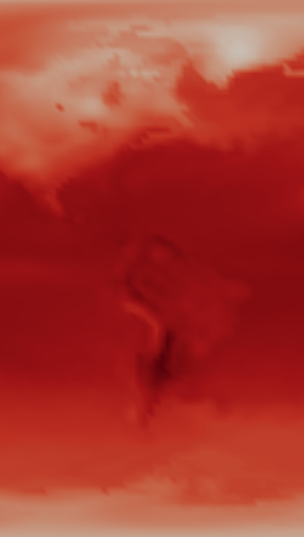}
  }
  \subfigureCaptionSkip
  \caption{%
    Climate dataset for $t=0$~(a), as well as some excerpts~(box) for subsequent timesteps. We can
    see that the continent of Africa exhibits a temperature increase for $t=3$ and $t=4$.
  }
  \label{fig:Climate data}
\end{figure}
\begin{figure}[tbp]
  \centering
  \subfigure[Wasserstein distance]{%
    \iffinal
      \includegraphics{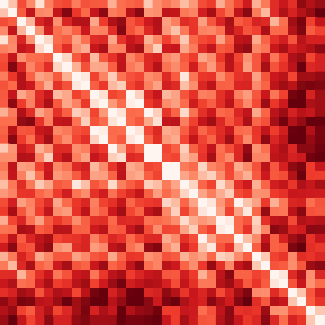}
    \else
      \begin{tikzpicture}
        \begin{axis}[%
            enlargelimits     = false,
            axis x line       = none,
            axis y line       = none,
            width             = 8cm,
            height            = 8cm,
            unit vector ratio =1 1 1,
          ]
          \addplot[%
            matrix plot,
            point meta     = explicit,
            mesh/ordering  = colwise,
            faceted color  = none,
            shader         = faceted,
            point meta min = 3.50,
            point meta max = 4.75,
            colormap/Reds,
            ]
            table[header=false, meta index=2] {Data/DKRZ/Distances_Wasserstein_2_transformed.txt};
          \end{axis}
      \end{tikzpicture}
    \fi
  }
  \quad
  \subfigure[\Gls{ISPH} distance]{%
    \iffinal
      \includegraphics{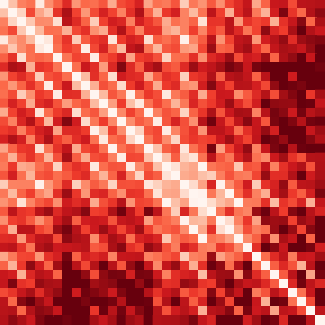}
    \else
      \begin{tikzpicture}
        \begin{axis}[%
            enlargelimits     = false,
            axis x line       = none,
            axis y line       = none,
            width             = 8cm,
            height            = 8cm,
            unit vector ratio =1 1 1,
          ]
          \addplot[%
            matrix plot,
            point meta     = explicit,
            mesh/ordering = colwise,
            faceted color  = none,
            shader         = faceted,
            point meta min = 400000,
            point meta max = 600000,
            colormap/Reds,
            ]
            table[header=false, meta index=2] {Data/DKRZ/Distances_extended_persistence_hierarchy_transformed.txt};
          \end{axis}
      \end{tikzpicture}
    \fi
  }
  \subfigureCaptionSkip
  \caption{%
    Dissimilarity matrices comparing the Wasserstein distance and our proposed \gls{ISPH} distance,
    which is capable of detecting the oscillatory behavior~(minor diagonals) inherent to the
    data.
  }
  \label{fig:Climate data dissimilarity matrices}
\end{figure}
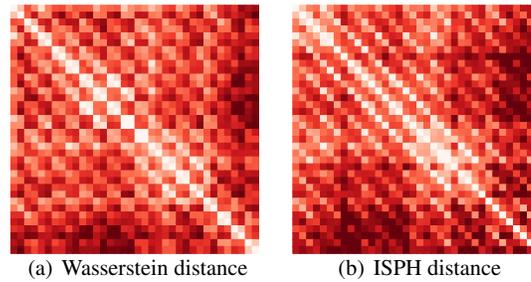

%
We used time-varying scalar field data~(surface temperature) from the German Climate Computing
Center~(DKRZ). The large size~(18,432 positions) and the large number of timesteps~(1460) makes
comparing these data complicated. Fig.~\ref{fig:Climate data} shows an excerpt of the dataset. It exhibits oscillatory behavior because~(global) temperature follows a day--night pattern. At the given resolution, a full day--night cycle comprises four timesteps. We would hence expect that the topological dissimilarity between corresponding timesteps is somewhat equal; or, more precisely, we would expect that a pairwise distance matrix depicts the oscillatory behavior that occurs in the data.
To assess the capabilities of the dissimilarity measure from Sec.~\ref{sec:Dissimilarity measure},
we compare it to the Wasserstein distance.
We first calculate pairwise distances between the first 36 timesteps of
the data. Calculating the Wasserstein distance takes about \SI{2}{\minute} per pair, hence comparing
all 36 pairs takes about \SI{21}{\hour}. Our hierarchy-based dissimilarity measure, by contrast,
takes \SI{2.3}{\second} for calculating each timestep and approximately \SI{6}{\second} to calculate
the dissimilarity per pair. The full distance matrix is thus obtained in approximately
\SI{1}{\hour}.
Fig.~\ref{fig:Climate data dissimilarity matrices} depicts the dissimilarity matrices of the first
36 timesteps. We can see that the \gls{ISPH} distance matrix contains patterns in the minor
   diagonals, indicating that timesteps $t_i$ and $t_{i+4}$ are highly similar. These patterns appear because the corresponding \glspl{ISPH} are also highly similar, even though the persistence pairs~(and thus the persistence diagrams) change over time because the range of temperatures changes. The Wasserstein distance does not exhibit these patterns. We note, however, that the two matrices are highly correlated~($R^2 \approx 0.95$), showing that for \emph{most} of the timesteps, both measures yield similar values.

\section{Conclusion and Future Work}

We presented a novel hierarchy that relates persistence pairs with each other. In contrast to
earlier work, our hierarchy is better capable of distinguishing certain nesting behaviors, which
we demonstrated by means of several example datasets.
At present, it is unclear to what extent the persistence hierarchy is an invariant like the
\emph{rank invariant} of multidimensional persistence~\cite{Carlsson09}.
Since the hierarchy changes when the data undergoes certain transformations, it remains to be shown
which operations leave it unchanged---a trivial observation is that the hierarchy is invariant with
respect to uniform scaling in all critical points.
Hence, the development of metrics for matching hierarchies may be beneficial. We only briefly
sketched a dissimilarity measure based on labeled trees. By considering the complete connectivity
structure of the tree, though, graph kernels such as the random walk kernel~\cite{Vishwanathan10}
could be employed.
Furthermore, it would be interesting to analyze how the hierarchy changes when different pairing
schemes for critical points are employed, such as the measures proposed by Carr et al.~\cite{Carr04}.

\begin{acknowledgement}
  We thank the anonymous reviewers for their helpful comments that helped us improve and clarify
  this manuscript.
\end{acknowledgement}


\bibliographystyle{spmpsci}
\bibliography{TopoInVis2017_Hierarchies}

\end{document}